\documentstyle[12pt,twoside]{article}
\input epsf.tex
\newcommand{\semi}{\mbox{$\times\!\rule{.2mm}{2mm}$}}
\newcommand{\NN}{\mbox{$I\!\! N$}}  

\let\\=\cr

\def\Chi{\hbox{\raise0.5ex\hbox{$\chi$}}}

\newtheorem{th}{Theorem}

\newtheorem{prop}{Proposition}

\newtheorem{defn}{Definition}

\newtheorem{rem}{Remark}

\input epsf.tex

\def\picill#1by#2(#3)
{ \begin{center}\leavevmode \epsfxsize = #1 \epsffile{#3}
\end{center}}

\begin{document}
\pagestyle{myheadings}
\markboth{{\sc S. Lambropoulou}}{{\sc Braid structures in 3--manifolds}}

\title{Braid structures in knot complements, handlebodies and 3--manifolds }

\author{Sofia Lambropoulou  \\
Mathematisches Institut, G\"ottingen Universit\"at}

\date{}

\maketitle

\begin{abstract}

We consider braids on $m+n$ strands, such that the first $m$ strands  are trivially fixed.
 We denote the set of all such braids by $B_{m,n}$. Via
concatenation $B_{m,n}$ acquires a group structure. The objective of this paper is to find a
presentation for  $B_{m,n}$ using  the structure of its corresponding pure braid subgroup,
$P_{m,n}$, and the fact that it is a subgroup of the classical Artin group $B_{m+n}$.  Then we
give an irredundant presentation for  $B_{m,n}$. The paper concludes by showing that these braid
groups or appropriate cosets of them are related to knots in handlebodies, in knot complements and in
c.c.o. 3--manifolds.

\end{abstract}

\section{Introductory notions and motivations}

\begin{defn}{ \ The set of all elements of the classical Artin group $B_{m+n}$ for which, if we
remove the last $n$ strands we are left with the identity braid on $m$ strands, shall be denoted by
$B_{m,n}$ (see figure 1(a) below for an example in  $B_{3,3}$). The elements of $B_{m,n}$ are special
cases of `mixed braids' (cf. section 6).} \end{defn}

\begin {center}
\leavevmode \epsfxsize =13cm \epsffile{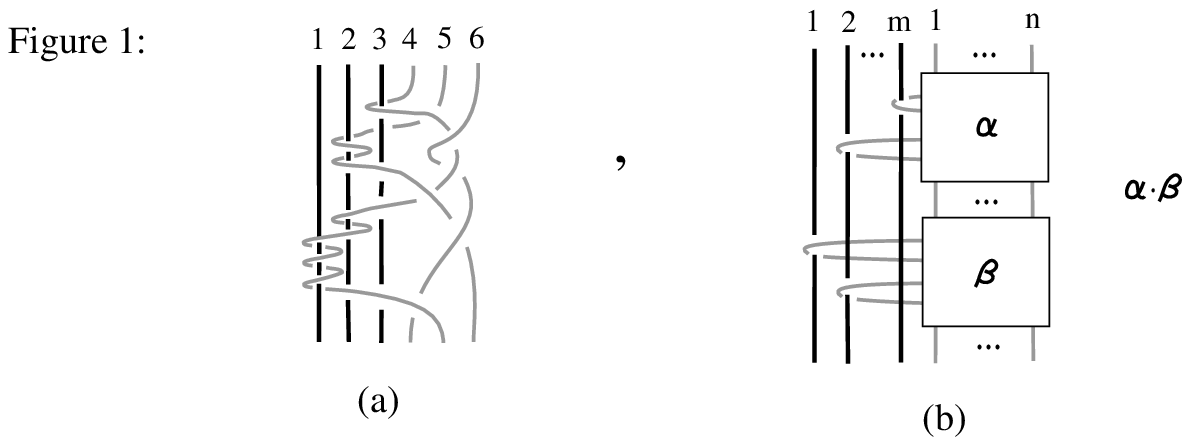}
\end {center}

\noindent Concatenation is a closed operation in $B_{m,n}$: the product  $\alpha\cdot
\beta$  of two elements   $\alpha, \beta\in B_{m,n}$ is also an element of $B_{m,n}$ (see
figure 1(b)).  Thus $B_{m,n} \leq B_{m+n}$. Our purpose is to obtain a simple presentation for
$B_{m,n}$.

\bigbreak

The motivation for studying these braids comes from studying
oriented knots and links in knot  complements, in c.c.o.
3-manifolds and in handlebodies, since  these spaces may be
represented by a fixed braid or a fixed integer-framed braid in
$S^3$. Then knots and links in these spaces may be represented by
elements of the above  braid groups $B_{m,n}$ or of appropriate
cosets of these groups. More precisely, if  $M$ denotes the
complement of the $m$-unlink or a connected sum of $m$ lens spaces
of  type $L(p,1)$ or a handlebody of genus $m$, then knots and
links in these spaces may be represented precisely by the  mixed
braids in $B_{m,n}$, for $n \in \NN$. In the case $m=1,$ $B_{1,n}$
is the Artin group of type $\cal B$ (cf. [4], [5], [6]). If $M$ is
generic, concatenation is no more a closed operation of mixed
braids, but as we show in section 6, knots and links in $M$ may be
represented  by   mixed braids in $B_{m,n}$, for $n \in \NN$,
followed by a fixed part associated to $M$, i.e. by elements of a
coset of $B_{m,n}$.

\bigbreak

We recall now some facts about braids and pure braids. For more
details and a complete study of the classical theory of braids the
reader is referred to [1].  The
 pure braid group, $P_n$, corresponding to the classical Artin group on $n$ strands, $B_n$,
consists of all elements in $B_n$ that induce the identity permutation in $S_n$, thus $P_n \lhd B_n$
and $P_n$ is generated by the elements

\vspace{.1in}
$ \begin{array}{ll}
a_{ij} & =  {\sigma_i}^{-1}{\sigma_{i+1}}^{-1}\ldots
{\sigma_{j-2}}^{-1}{\sigma_{j-1}}^2 \sigma_{j-2}\ldots\sigma_{i+1}\sigma_i \\ [.1in]
  & =  \sigma_{j-1}\sigma_{j-2}\ldots \sigma_{i+1}{\sigma_i}^2
{\sigma_{i+1}}^{-1}\ldots {\sigma_{j-2}}^{-1}{\sigma_{j-1}}^{-1}, \ \  1\leq i<j\leq n.
\end{array} $
\vspace{.1in}

\noindent The generators $a_{ij}$ may be pictured geometrically as an elementary loop between the
 $i$th and $j$th strand (cf. figure 2).

\bigbreak

\noindent The most important property of pure braids is that they have a canonical form, the
so-called  {\it `Artin's canonical form'}, which says that every element, $A$, of $P_n$ may be
written uniquely in the form:
$$ A=U_1U_2 \cdots U_{n-1} $$
where each $U_i$ is a uniquely determined product of powers of the $a_{ij}$ using
only those with $i<j$. Geometrically, this means that any pure braid can be `combed' i.e. can
be written canonically as the pure braiding of the first strand with the rest, then keep the
first  strand fixed and uncrossed and have the pure braiding of the second strand and so on (cf.
figure 3).

\vspace{.1in}

The main idea for finding a presentation for $P_n$ is the following: The combing of a strand may
be regarded as a loop in the complement space of the other strands, and as such is an element of a
free group  since the fundamental group of a punctured disc is free. Thus,
 $$ P_n =  F_{n-1} \semi \cdots \semi F_2 \semi F_1 = F_{n-1} \semi P_{n-1}, $$
where each $F_i$ is a free group on the generators \ $a_{1,i+1}, \ldots, a_{i,i+1}$ \
(the elementary loops between the $(i+1)$st strand and all its previous ones), and where the
action is induced by conjugation. It turns out that $P_n$ has $\frac{n(n-1)}{2}$ generators and
 \ $\frac{1\cdot 2^2 + 2\cdot 3^2 + \ldots + (n-2)(n-1)^2}{2}$ \ relations, the following:

\[ {a_{ij}}^{-1}{a_{rs}}{a_{ij}}= \left\{ \begin{array}{ll}
{a_{rs}} \ \ \ \ \ \ \ \ \ \ \ \ \ \ \ \ \ \ \ \ \ \ \  \mbox{ \ if \ $i<j<r<s$ } &
\mbox{ \ or \ } r<i<j<s, \\
{a_{is}}{a_{js}}{a_{is}}^{-1} & \mbox{ \ if \ } i<j<s, \\
{a_{is}}{a_{js}}{a_{is}}{a_{js}}^{-1}{a_{is}}^{-1} & \mbox{ \ if \ } i<j<s, \\
{a_{is}}{a_{js}}{a_{is}}^{-1}{a_{js}}^{-1}{a_{rs}}{a_{js}}{a_{is}}{a_{js}}^{-1}
{a_{is}}^{-1} & \mbox{ \ if \ } i<r<j<s.  \end{array} \right. \]

\vspace{.1in}


Based on these ideas, we introduce in section 2 the pure braid group $P_{m,n}$ and in section 3 we
find a presentation for it. Then in section 4 we put together a presentation for $B_{m,n}$,
which
  we simplify in section 5. In section 5 we also give a Dynkin-diagram related to $B_{m,n}$. Finally,
in section 6 we explain that elements of $B_{m,n}$ represent oriented knots and links in certain
spaces and that appropriate cosets of $B_{m,n}$ represent knots and links in the generic cases of
knot complements and c.c.o. 3-manifolds.

\bigbreak

\noindent The results here have been preliminary studied by the
author in [4] and have been presented in various mathematical
meetings since 1995. A. Sossinsky, independently, motivated
 by the same topological considerations,  studies these groups in [8] and he conjectures the
irredundant presentation for $B_{m,n}$. Moreover, V. Vershinin in
[9] studies the  groups $B_{m,n}$ in connection to handlebodies of
genus $g$, taking a configuration-spaces approach. Back in 1993
Alastair Leeves had found a presentation for $B_{m,n}$, which was
presented in [4], but a proof was never published. The author is
thankfull to A. Leeves for inspiring discussions  at the time.
Also, her grateful thanks are due to Bernard Leclerc for his
careful reading through this work and his very valuable comments.

\bigbreak

 \noindent This is the first paper in a sequel of three. The next one gives  expressions for
 algebraic equivalence of braids reflecting knot isotopy in arbitrary knot complements and c.c.o.
 3-manifolds. The case of handlebodies is joint work with Reinhard H\"aring-Oldenburg.

\section{  The pure braid group $P_{m,n}$}

\begin{defn}{ \ The corresponding pure braid group $P_{m,n}$ of $B_{m,n}$ is defined as
$P_{m,n} = B_{m,n} \cap P_{m+n}$, i.e. $P_{m,n} \leq P_{m+n}$ and it does not contain pure braiding
among the first $m$  strands. } \end{defn}

\noindent By its definition, $P_{m,n}$ is generated by the pure braid generators $a_{ij}$ for $i
\in \{1, \ldots, m+n-1\}$ and  $j \in \{m+1, \ldots, m+n\}$ of $P_{m+n}$ (see figure 2).
Then, $B_{m,n}$ is clearly generated by the elementary mixed braids (drawn below) $a_{ij}$ for $i
\in \{1, \ldots, m+n-1\}$ and $j \in \{m+1, \ldots, m+n\}$ together with $\sigma_{m+1}, \ldots,
\sigma_{m+n-1}$, the elementary crossings among the last $n$ strands. Note that the inverses of
the $a_{ij}$'s and $\sigma_k$'s  are represented by the same geometric  pictures, but with the
opposite crossings.

\begin {center}
\leavevmode \epsfxsize =13cm \epsffile{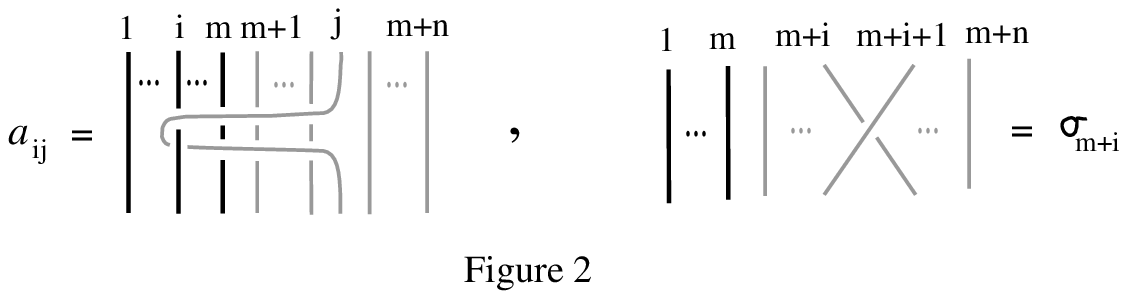}
\end {center}
Also, by definition we have an exact sequence

\[ 1 \longrightarrow P_{m,n} \longrightarrow B_{m,n} \longrightarrow S_n \longrightarrow 1.\]

In particular, $P_{m,n} \lhd B_{m,n}$. More precisely, we have the following relations:

\[ {\sigma_k}^{-1}{a_{ij}}^{\pm}{\sigma_k} = \left\{ \begin{array}{ll}
{a_{ij}}^{\pm}  \mbox{ \ \ \ \ \ if \ $k\leq i-2$} & \mbox{ or \ $i+1
\leq k \leq j-2$ \ or \  $k \geq j+1,$ } \\
{a_{i-1,j}}^{\pm} &  \mbox{  \ if \ $k=i-1,$} \\
{a_{ij}}{a_{i+1,j}}^{\pm}{a_{ij}}^{-1}  &  \mbox{ \ if \ $k=i$,} \\
{a_{i,j-1}}^{\pm} & \mbox{  \ if \ $k=j-1$,} \\
{a_{ij}}{a_{i,j+1}}^{\pm}{a_{ij}}^{-1}  &  \mbox{ \ if \ $k=j$. } \hfill \Box \end{array}
\right. \]

\noindent We shall call these  {\it mixed relations} and we shall denote them by $M_1,
M_2, M_3, M_4$ and $M_5$ in the order they are written.

\bigbreak

\noindent Thus $B_{m,n}$ is a group extension of $P_{m,n}$ by $S_n$. This  will yield a presentation
for $B_{m,n}$, conditionally to knowing a presentation for  $P_{m,n}$.

\section{A presentation for $P_{m,n}$}

\begin{th}{ \ The pure braid group $P_{m,n}$ is generated by the elements $a_{ij}$ for $i
\in \{1, \ldots, m+n-1\}$, $j \in \{m+1, \ldots, m+n\}$ and $i<j$, which are subject to the
relations:
\[\begin{array}{crclll}
(P_1) & {a_{ij}}^{-1}{a_{rs}}{a_{ij}}   & = & {a_{rs}} \mbox{ \ \ \ \ \ \ \ \ \ \ \ \ \ \ \ \ if \
} i<j<r<s  & \mbox{ or } & r<i<j<s,\\
(P_2) & {a_{ij}}^{-1}{a_{js}}{a_{ij}}   & = & {a_{is}}{a_{js}}{a_{is}}^{-1}    & \mbox{ if } &
i<j<s, \\
(P_3) & {a_{ij}}^{-1}{a_{is}}{a_{ij}}   & = &
{a_{is}}{a_{js}}{a_{is}}{a_{js}}^{-1}{a_{is}}^{-1}  & \mbox{ if } & i<j<s, \\
(P_4) & {a_{ij}}^{-1}{a_{rs}}{a_{ij}}   & = &
{a_{is}}{a_{js}}{a_{is}}^{-1}{a_{js}}^{-1}{a_{rs}}{a_{js}}  {a_{is}}{a_{js}}^{-1}{a_{is}}^{-1}  &
\mbox{ if } & i<r<j<s. \end{array}\]
} \end{th}
Relations $P_1, P_2, P_3$ and $P_4$ shall be called {\it pure braid relations}. Note that
$P_1$ and $P_4$ involve the strands $i, r, j, s$, whilst $P_2$ and $P_3$ involve the strands
$i, j, s$.
\bigbreak
\noindent {\bf Proof} \ By its definition and by the fact that the $a_{ij}$'s (for all indices)
generate $P_{m+n}$ follows that the above set of elements is indeed a set of generators for
$P_{m,n}$. Since $P_{m,n} \leq P_{m+n}$ we can apply on its elements Artin's combing. As for $P_n$,
the combing of a strand can be regarded as a loop in the complement space of the strands with
smaller index (including the $m$ fixed ones) and as such it is an element of a free group.
Therefore we have:

\vspace{.1in}
\noindent $ P_{m,1} = F_m = \langle a_{1,m+1}, a_{2,m+1}, \ldots, a_{m,m+1}
\rangle,$ the free group on $m$ generators.

\vspace{.1in}
\noindent  Further is:
$ P_{m,2} = F_{m+1} \semi P_{m,1} = \langle a_{1,m+2}, \ldots, a_{m,m+2},
a_{m+1,m+2} \rangle \semi P_{m,1}, $ \ i.e. $F_{m+1}$ is the free group on $m+1$ generators, and
$P_{m,1}$ acts on $F_{m+1}$ by conjugation,  via the relations of the pure braid group $P_{m+1}$
for appropriate indices.  We proceed inductively to obtain:

\[\begin{array}{rcl}
P_{m,n}   & = & F_{m+n-1} \semi \cdots \semi F_{m+1} \semi F_m  = F_{m+n-1} \semi P_{m,n-1} \\
          & = & \langle a_{1,m+n}, \ldots, a_{m,m+n},
                \ldots, a_{m+n-1,m+n} \rangle \semi P_{m,n-1},  \end{array}\]

\noindent where  $F_{m+n-1}$ is the free group on $m+n-1$
generators, and where  $P_{m,n-1}$ acts on $F_{m+n-1}$ by
conjugation,  via the relations of the pure braid group $P_{m+n}$
for appropriate indices, i.e. via the relations  $P_1, P_2, P_3$
and $P_4$. $\hfill \Box$

\bigbreak
Some remarks are now due.

\begin{rem}{\rm \ The groups $P_{m,n}$ and $P_{m+n}$ have seemingly the same presentation. For
$m\neq 1$ is, though, $P_{m,n} \neq P_{m+n}$. The difference lies in the restriction of the
indices of the generators. In fact, $P_{m,n}$ has  $\frac{n(n+2m-1)}{2}$ generators, which is the
number of generators of $P_{m+n}$ less the number of generators of $P_m$. Moreover, $P_{m,n}$ has
$ \ \frac{(m-1)\cdot m^2 + m\cdot (m+1)^2 + \cdots + (m+n-2)(m+n-1)^2}{2} - \frac{(m-1)\cdot
m\cdot n\cdot (n+2m-1)}{4}$ \ relations, which is the number of relations of $P_{m+n}$ less the
number of relations of $P_m$. In the case $m = 1$  holds  $P_{1,n} = P_{1+n},$ which
follows immediately  from the definition of $P_{m,n}$ or can be observed from its presentation for
$m = 1$. } \end{rem}

\begin{rem}{\rm \ In [4] there is a discussion about the groups $B_{m,n}$ and
a different line of proof is given for finding a presentation.
There, by $P_{m,n}$ we denoted some smaller pure braid subgroups,
for which it is rather complicated to find a presentation. But the
case $m = 1$ was extensively treated, also in the sequel papers
[5, 6]. In all these previous results  $P_{1,n}$ denoted the free
group  $ F_n = \langle a_{12}, a_{13}, \ldots, a_{1,n+1} \rangle$
and {\it not} the corresponding pure braid group of $B_{1,n}$.
That's why we had then $B_{1,n}= F_n \semi B_n $. We hope that the
readers familiar with those results will not be in confusion.
 } \end{rem}

\section{A presentation for $B_{m,n}$}

In section 2 we showed that $ 1 \longrightarrow P_{m,n} \longrightarrow B_{m,n} \longrightarrow
S_n \longrightarrow 1$ and in section 3 we found a presentation for $P_{m,n}$.  Recall that
$S_n$ has the  presentation:

\[ \langle s_1,...,s_{n-1} \ | \ s_i s_j = s_j s_i  \mbox{ \ for \ } |i-j|>1,\ \ s_i s_{i+1}
           s_i = s_{i+1} s_i s_{i+1},\ \ {s_i}^2 = 1     \rangle .\]

\noindent We are now ready to put together a presentation for
$B_{m,n}$. Namely, we can apply a result from the theory of group
presentations (see [3], p.139), that gives a presentation for a
group extension  of two groups with known presentations. Indeed,
the following is then a presentation for $B_{m,n}$.

\[  \left< \begin{array}{ll}  \begin{array}{l}
a_{1,m+1}, \ldots, a_{1,m+n}, \ldots, a_{m,m+1}, \ldots, a_{m,m+n}, \\
a_{m+1,m+2}, \ldots, a_{m+1,m+n}, \ldots, a_{m+n-1,m+n}, \\
\sigma_{m+1}, \sigma_{m+2}, \ldots ,\sigma_{m+n-1} \\
\end{array} &
\left|
\begin{array}{l} P_1, P_2, P_3, P_4,       \\
                 M_1, M_2, M_3, M_4, M_5,  \\
\Sigma_1, \Sigma_2, \Sigma_3.              \\
\end{array}
\right.  \end{array} \right>,  \]

\noindent where the relations $\Sigma_1, \Sigma_2$ and $\Sigma_3$ are satisfied by the
$\sigma_{m+1}, \sigma_{m+2}, \ldots ,\sigma_{m+n-1}$ and they are the following:

\[\begin{array}{lrclll}
(\Sigma_1) \ \ \ \ & \sigma_i\sigma_j             & = & \sigma_j\sigma_i & \mbox{ if  } & |i-j|>1,
\\ (\Sigma_2) \ \ \ \ & \sigma_i\sigma_{i+1}\sigma_i & = & \sigma_{i+1}\sigma_i\sigma_{i+1} &
\mbox{ if  } &   m+1 \leq i \leq m+n-2, \\
(\Sigma_3) \ \ \ \ & {\sigma_i}^2      & = & a_{i,i+1} & \mbox{ if  } &  m+1 \leq i \leq m+n-2.
 \end{array}\]

\noindent  $\Sigma_1$ and $\Sigma_2$ are the  {\it `braid relations'}.

\bigbreak

\noindent Notice now that  relations $\Sigma_3$ for  $i \in \{ m+1, \ldots, m+n-1 \}$ and $j
\in \{ m+2, \ldots, m+n \}$ do not involve any mixed braiding and so they may be taken as defining
relations, namely:
$$\begin{array}{lcl}
{a_{m+1,m+2}}^{\pm} & := & {\sigma_{m+1}}^{\pm 2},  \\
{a_{m+1,m+3}}^{\pm} & := & \sigma_{m+2} {\sigma_{m+1}}^{\pm 2} {\sigma_{m+2}}^{-1}, \ldots  \\

{a_{ij}}^{\pm} & := & \sigma_{j-1}\ldots {\sigma_{i+1}}{\sigma_i}^{\pm 2} {\sigma_{i+1}}^{-1}\ldots
{\sigma_{j-1}}^{-1},   \\

\ &  \vdots & \     \\
{a_{m+n-1,m+n}}^{\pm} & := & {\sigma_{m+n-1}}^{\pm 2}.   \\
\end{array}$$
\noindent Therefore, we want to omit eventually these $a_{ij}$'s
 from the list of generators of $B_{m,n}$ and subsequently to eliminate or simplify
all relations involving these elements, applying Tietze
transformations. Indeed, we examine one by one the relations and we have:

\vspace{.1in}
\noindent $\mbox{\boldmath $P_1:$}  \ \ {a_{ij}}^{-1}{a_{rs}}{a_{ij}}    =  {a_{rs}}$ for the
case $ r<i<j<s$. If all  $r, i, j, s \in \{m+1, \ldots, m+n\}$, the relations follow from
$\Sigma_1$ and $ \Sigma_2$  and so we only keep the ones where  $r \in \{1, \ldots, m\}$ and  $i
\in \{m+1, \ldots, m+n-1\}$ or $r, i \in \{1, \ldots, m\}$.

\vspace{.1in}
\noindent $\mbox{\boldmath $P_2:$}  \ \ {a_{ij}}^{-1}{a_{js}}{a_{ij}} =
{a_{is}}{a_{js}}{a_{is}}^{-1}$ for $i<j<s$.  Since  $j, s \in \{m+1, \ldots, m+n\}$ the only case
to be kept is when $i \in \{1, \ldots, m\}$, as the relations for $i \in \{m+1, \ldots, m+n-1\}$
follow from the braid relations.

\vspace{.1in}
\noindent $\mbox{\boldmath $P_3:$}  \ \ {a_{ij}}^{-1}{a_{is}}{a_{ij}}  =
{a_{is}}{a_{js}}{a_{is}}{a_{js}}^{-1}{a_{is}}^{-1}$ for $i<j<s.$ As in the previous case, the
only relations that do not follow from $\Sigma_1$ and $ \Sigma_2$ are the ones where $i \in \{1,
\ldots, m\}$.

\vspace{.1in}
\noindent $\mbox{\boldmath $P_4:$}  \ \ {a_{ij}}^{-1}{a_{rs}}{a_{ij}} =
{a_{is}}{a_{js}}{a_{is}}^{-1}{a_{js}}^{-1}{a_{rs}}{a_{js}}  {a_{is}}{a_{js}}^{-1}{a_{is}}^{-1}$
for $i<r<j<s$. Since  $j, s \in \{m+1, \ldots, m+n\}$ the only relations to be kept are those
where either $i \in \{1, \ldots, m\}$ and $r \in \{m+1, \ldots, m+n-1\}$ or $i, r \in \{1, \ldots,
m\}$.

\vspace{.1in}
\noindent $\mbox{\boldmath $M_1:$}  \ \ {\sigma_k}^{-1}{a_{ij}}^{\pm}{\sigma_k} = {a_{ij}}^{\pm}$
for $k\leq i-2$ or
 $i+1 \leq k \leq j-2$  or $k \geq j+1$.  Here also we have $j \in \{m+1, \ldots, m+n\}$ and
$k \in \{m+1, \ldots, m+n-1\}$. Now, if $i \in \{m+1, \ldots, m+n-1\}$, all these relations
follow from $\Sigma_1$ and $\Sigma_2$, whilst for $i \in \{1, \ldots, m\}$ it only makes sense to
consider  $k\leq j-2$ or $k \geq j+1$.

\vspace{.1in}
\noindent $\mbox{\boldmath $M_2:$}  \ \ {\sigma_{i-1}}^{-1}{a_{ij}}^{\pm}{\sigma_{i-1}} =
{a_{i-1,j}}^{\pm}$. Since  $i-1 \in \{m+1, \ldots, m+n-1\}$ it must be $i > m+1$ and so all these
relations follow from  $\Sigma_1$ and $\Sigma_2$.

\vspace{.1in}
\noindent $\mbox{\boldmath $M_3:$}  \ \ {\sigma_i}^{-1}{a_{ij}}^{\pm}{\sigma_i} =
{a_{ij}}{a_{i+1,j}}^{\pm}{a_{ij}}^{-1}$. This is analogous to the above case, since  $i \geq m+1$.
\vspace{.1in}

\noindent $\mbox{\boldmath $M_4:$}  \ \ {\sigma_{j-1}}^{-1}{a_{ij}}^{\pm}{\sigma_{j-1}} =
{a_{i,j-1}}^{\pm}$. Here also the only cases that do not follow from the braid relations are the
ones with $i \in \{1, \ldots, m\}$.

\vspace{.1in}
\noindent $\mbox{\boldmath $M_5:$}  \ \ {\sigma_j}^{-1}{a_{ij}}^{\pm}{\sigma_j} =
{a_{ij}}{a_{i,j+1}}^{\pm}{a_{ij}}^{-1}$. As above, the only relations surviving
are the ones where $i \in \{1, \ldots, m\}$.

\begin{rem} {\rm \ The remaining relations of $P_1$ for $i \in \{m+1, \ldots, m+n-1\}$  follow from
the simpler relations: \ \ \ $ {a_{ij}} \sigma_k  = \sigma_k {a_{ij}}$ for $k\leq j-2$ or $k \geq
j+1$, which  coincides with the remaining of  $M_1$ above.  } \end{rem}

\noindent To summarize, we showed that the following is a presentation for $B_{m,n}$.

\[  B_{m,n} = \left< \begin{array}{ll}  \begin{array}{l}
a_{1,m+1}, \ldots, a_{1,m+n}, \ldots, a_{m,m+1}, \ldots, a_{m,m+n}, \\
a_{m+1,m+2}, \ldots, a_{m+1,m+n}, \ldots, a_{m+n-1,m+n}, \\
\sigma_{m+1}, \sigma_{m+2}, \ldots ,\sigma_{m+n-1} \\
\end{array} &
\left|
\begin{array}{l} P'_1, P'_2, P'_3, P'_4,       \\
                 M'_1, M'_2, M'_3,   \\
\Sigma_1, \Sigma_2              \\
\end{array}
\right.  \end{array} \right>,  \]

\noindent where we have:

\[\begin{array}{crcl}
(P'_1) \ \ & {a_{ij}} {a_{rs}}   & = & {a_{rs}} {a_{ij}} \ \ \ \ \ \ \ \ \ \ \ \ \ \ \ \ \ \mbox{
\ for \ }  r<i<j<s, \ 1 \leq i,r \leq m,   \\
(P'_2) \ \ & {a_{ij}}^{-1}{a_{js}}{a_{ij}}   & = &
{a_{is}}{a_{js}}{a_{is}}^{-1} \ \ \ \ \ \ \ \ \ \ \ \mbox{ \ for \ }  i<j<s, \ 1 \leq i \leq
m,   \\
(P'_3) \ \ & {a_{ij}}^{-1}{a_{is}}{a_{ij}}   & = &
{a_{is}}{a_{js}}{a_{is}}{a_{js}}^{-1}{a_{is}}^{-1} \ \mbox{ \ for \ }
i<j<s, \ 1 \leq i \leq m,    \\
(P'_4) \ \ & {a_{ij}}^{-1}{a_{rs}}{a_{ij}}   & = &
{a_{is}}{a_{js}}{a_{is}}^{-1}{a_{js}}^{-1}{a_{rs}}{a_{js}} {a_{is}}{a_{js}}^{-1}{a_{is}}^{-1}
 \mbox{ \ for \ }  i<r<j<s,  \\
 \ & \ & \ & 1 \leq i \leq m, \ 1 \leq r \leq m+n-1,    \\
(M'_1) \ \ & {\sigma_k}^{-1}{a_{ij}}^{\pm}{\sigma_k}   & = &  {a_{ij}}^{\pm} \ \
 \mbox{ \ for \ } k\leq j-2  \mbox{ \ or \ }  k \geq j+1 \mbox{ \ and \ } 1 \leq i \leq m,  \\
(M'_2) \ \ & {a_{ij}}^{\pm}  & = & {\sigma_{j-1}}{a_{i,j-1}}^{\pm}{\sigma_{j-1}}^{-1} \ \ \mbox{ \
for \ }  1 \leq i \leq m,     \\
(M'_3) \ \ & {\sigma_j}^{-1}{a_{ij}}^{\pm}{\sigma_j}   & =
& {a_{ij}}{a_{i,j+1}}^{\pm}{a_{ij}}^{-1} \ \ \mbox{ \ for \ }  1 \leq i \leq m.
 \end{array}\]

Having now done the first, `obvious' clearing in the original presentation of $B_{m,n}$, we
observe that many of the above relations are redundant or they simplify further, and that we may
omit the $a_{ij}$'s with $i \geq m+1$.

\begin{th}{ \ The following is a presentation for $B_{m,n}$:

\[  B_{m,n} = \left< \begin{array}{ll}  \begin{array}{l}
a_{1,m+1}, \ldots, a_{1,m+n}, \cdots,  \\
a_{m,m+1}, \ldots, a_{m,m+n}, \\
\sigma_{m+1}, \sigma_{m+2}, \ldots ,\sigma_{m+n-1} \\
\end{array} &
\left|
\begin{array}{l} \Sigma_1, \Sigma_2, (1), (2), (3), (4),       \\
                 \mbox{for all appropriate indices}   \\
\end{array}
\right.  \end{array} \right>,  \]

\noindent where we have:

\[\begin{array}{lcccl}
(1) \ \  & {\sigma_k}^{-1}{a_{ij}}^{\pm}{\sigma_k}   & = & {a_{ij}}^{\pm} &  \mbox{ for
\ } k\leq j-2  \mbox{ \ or \ }  k \geq j+1,  \\

(2) \ \  & {a_{ij}}^{\pm}    & = &   {\sigma_{j-1}}{a_{i,j-1}}^{\pm}{\sigma_{j-1}}^{-1}, & \ \\

(3) \ \  & {\sigma_j}^{-1}{a_{ij}}^{\pm}{\sigma_j}    & = &
{a_{ij}}{a_{i,j+1}}^{\pm}{a_{ij}}^{-1} & \    \\

(4) \ \  & {a_{ij}}^{\pm} {a_{r,j+1}}^{\pm}   & = &
{a_{r,j+1}}^{\pm} {a_{ij}}^{\pm}  & \mbox{ for \ } r<i.  \\

 \end{array}\]
 } \end{th}

\bigbreak
\noindent {\bf Proof} \ Relations $(1), (2)$ and $(3)$ are precisely  $M'_1, M'_2$ and
$M'_3$, whilst relations $(4)$ are a special case of relations $P'_1$. So we have to show that
$P'_2, P'_3, P'_4$ and $P'_5$ as well as the rest cases of $P'_1$ follow from
$\Sigma_1, \Sigma_2, (1), (2), (3)$ and $(4)$. Before continuing we note that, using $(2)$,
relation $(3)$ is equivalent to
$$ {\sigma_j}{a_{ij}}{\sigma_j}{a_{ij}}^{\pm}  =
 {a_{ij}}^{\pm}{\sigma_j}{a_{ij}}{\sigma_j} $$
and relation $(4)$ is equivalent to
$${a_{ij}}^{\pm} (\sigma_j {a_{rj}}^{\pm}{\sigma_j}^{-1})   =
 (\sigma_j {a_{rj}}^{\pm}{\sigma_j}^{-1}) {a_{ij}}^{\pm}.$$
We shall also use  these forms in the proof. We proceed now case by
case. The underlining indicates the expressions involved in each step of the proof.

\vspace{.1in}

\noindent $\begin{array}{rcl}
 {a_{ij}} \underline{ {a_{rs}}} & \stackrel{M'_2}{=}  &

\underline{ {a_{ij}} \sigma_{s-1} \ldots \sigma_{j+1} } \underline{ {a_{r,j+1}} }
{\sigma_{j+1}}^{-1} \ldots {\sigma_{s-1}}^{-1} \\

 \ & \stackrel{M'_1}{=} & \sigma_{s-1} \ldots \sigma_{j+1} \underline{ {a_{ij}} (\sigma_j
{a_{rj}}{\sigma_j}^{-1}) } {\sigma_{j+1}}^{-1} \ldots {\sigma_{s-1}}^{-1} \\

 \ & \stackrel{(4)}{=} & \sigma_{s-1} \ldots \sigma_{j+1}  (\sigma_j {a_{rj}}{\sigma_j}^{-1})
{a_{ij}}  \underline{ {\sigma_{j+1}}^{-1} \ldots {\sigma_{s-1}}^{-1}} \\

 \ & \stackrel{M'_1}{=} & \underline{  \sigma_{s-1} \ldots \sigma_j
{a_{rj}}{\sigma_j}^{-1} \ldots {\sigma_{s-1}}^{-1} } {a_{ij}}  \\

 \ & \stackrel{M'_2}{=} & {a_{rs}} {a_{ij}}.  \\
\end{array}$

\vspace{.1in}

\noindent $\begin{array}{rcl}
 {a_{ij}}^{-1} & {a_{js}} & {a_{ij}} \ \stackrel{M'_2}{=}  \

{a_{ij}}^{-1} ( \underline{ \sigma_{s-1} \ldots \sigma_{j+1} }{\sigma_j}^2  \underline{
{\sigma_{j+1}}^{-1} \ldots {\sigma_{s-1}}^{-1} }) {a_{ij}}     \\

 \ & \stackrel{M'_1}{=} & ( \sigma_{s-1} \ldots \sigma_{j+1})  \underline{  {a_{ij}}^{-1}
{\sigma_j}^2 } {a_{ij}} ({\sigma_{j+1}}^{-1} \ldots {\sigma_{s-1}}^{-1})  \\

 \ & \stackrel{M'_3}{=} & ( \sigma_{s-1} \ldots \sigma_{j+1}) \sigma_j {a_{ij}} {a_{i,j+1}}^{-1}
 \underline{ {a_{ij}}^{-1} \sigma_j } {a_{ij}} ({\sigma_{j+1}}^{-1} \ldots {\sigma_{s-1}}^{-1})
\\

 \ & \stackrel{M'_3}{=} & ( \sigma_{s-1} \ldots \sigma_{j+1}) \sigma_j {a_{ij}}  \underline{
{a_{i,j+1}}^{-1}  \sigma_j } {a_{ij}} {a_{i,j+1}}^{-1} {a_{ij}}^{-1} {a_{ij}} ({\sigma_{j+1}}^{-1}
\ldots {\sigma_{s-1}}^{-1})  \\

 \ & \stackrel{M'_2}{=} & ( \sigma_{s-1} \ldots \sigma_{j+1}) \sigma_j {a_{ij}} {\sigma_j}^2
 {a_{ij}}^{-1} {\sigma_j}^{-1} ({\sigma_{j+1}}^{-1} \ldots {\sigma_{s-1}}^{-1})  \\ [1.8mm]

 \ & =  &  \underline{ ( \sigma_{s-1} \ldots \sigma_j) {a_{ij}} ({\sigma_j}^{-1}
\ldots {\sigma_{s-1}}}^{-1} \underline{ \sigma_{s-1} \ldots \sigma_j) {\sigma_j}^2
({\sigma_j}^{-1} \ldots {\sigma_{s-1}}^{-1} }  \\  [1.8mm]

 \ & \ & \times \, \underline{ \sigma_{s-1} \ldots \sigma_j)
 {a_{ij}}^{-1} ({\sigma_j}^{-1}  \ldots {\sigma_{s-1}}^{-1}) }      \\

 \ & \stackrel{M'_2}{=} & {a_{is}}{a_{js}}{a_{is}}^{-1}.   \\
\end{array}$

\vspace{.1in}

\noindent For $P'_3$ we have:

\vspace{.1in}

\noindent $\begin{array}{rcl}
 {a_{is}} &  \underline{ {a_{js}} } & {a_{is}} \ \underline{ {a_{js}}^{-1} } \ {a_{is}}^{-1}

\ \stackrel{M'_2}{=} \ \underline{ {a_{is}} } (\sigma_{s-1} \ldots \sigma_{j+1} {\sigma_j}^2
{\sigma_{j+1}}^{-1} \ldots {\sigma_{s-1}}^{-1})  \underline{ {a_{is}} } \\  [1.8mm]

\ & \ & \times \, (\sigma_{s-1} \ldots
\sigma_{j+1} {\sigma_j}^{-2}  {\sigma_{j+1}}^{-1} \ldots {\sigma_{s-1}}^{-1})  \underline{
{a_{is}}^{-1} } \\

 \ & \stackrel{M'_2}{=} &  (\sigma_{s-1} \ldots \sigma_j {a_{ij}}
{\sigma_j}^{-1} \ldots {\sigma_{s-1}}^{-1})  (\sigma_{s-1} \ldots \sigma_{j+1}
{\sigma_j}^2  {\sigma_{j+1}}^{-1} \ldots {\sigma_{s-1}}^{-1})  \\  [1.8mm]

\ & \ & \times \, (\sigma_{s-1} \ldots \sigma_j
{a_{ij}}  {\sigma_j}^{-1} \ldots {\sigma_{s-1}}^{-1})
 (\sigma_{s-1} \ldots \sigma_{j+1} {\sigma_j}^{-2}  {\sigma_{j+1}}^{-1} \ldots
{\sigma_{s-1}}^{-1})   \\   [1.8mm]

 \ & \ & \times \, (\sigma_{s-1} \ldots \sigma_j {a_{ij}}^{-1}
{\sigma_j}^{-1} \ldots {\sigma_{s-1}}^{-1})  \\   [1.8mm]

 \ & = &  (\sigma_{s-1} \ldots \sigma_j) {a_{ij}} {\sigma_j}^2 \underline{ {a_{ij}}
{\sigma_j}^{-2}  {a_{ij}}^{-1} } ({\sigma_j}^{-1} \ldots {\sigma_{s-1}}^{-1}).  \\
\end{array}$

\vspace{.1in}

\noindent On the other hand:

\vspace{.1in}

\noindent $\begin{array}{rcl}
 {a_{ij}}^{-1} & \underline{ {a_{is}} } & {a_{ij}} \ \stackrel{M'_2}{=} \
\underline{ {a_{ij}}^{-1} \sigma_{s-1} \ldots \sigma_{j+1} } \sigma_j {a_{ij}}
{\sigma_j}^{-1} \underline{ {\sigma_{j+1}}^{-1} \ldots {\sigma_{s-1}}^{-1} {a_{ij}} }  \\

 \ & \stackrel{M'_1}{=} &
 \sigma_{s-1} \ldots \sigma_{j+1} \underline{ {a_{ij}}^{-1} \sigma_j {a_{ij}} {\sigma_j}^{-1} }
{a_{ij}}  {\sigma_{j+1}}^{-1}  \ldots {\sigma_{s-1}}^{-1}  \\

 \ & \stackrel{M'_3}{=} &
 \sigma_{s-1} \ldots \sigma_{j+1} \sigma_j {a_{ij}} \sigma_j (\sigma_j {\sigma_j}^{-1})
{a_{ij}}^{-1} {\sigma_j}^{-2} {a_{ij}} ( \sigma_j {\sigma_j}^{-1}) {\sigma_{j+1}}^{-1}  \ldots
{\sigma_{s-1}}^{-1}  \\  [1.8mm]

 \ & = &
 (\sigma_{s-1} \ldots \sigma_j) {a_{ij}} {\sigma_j}^2
\underline{ {\sigma_j}^{-1} {a_{ij}}^{-1} {\sigma_j}^{-2} {a_{ij}} \sigma_j } ({\sigma_j}^{-1}
\ldots {\sigma_{s-1}}^{-1}).  \\
\end{array}$

\vspace{.1in}

\noindent Therefore, it suffices to show that the underlined expressions in the last equation of
either side are equal. Indeed we have:

\vspace{.1in}

\noindent $\begin{array}{rcl}
  {a_{ij}} & {\sigma_j}^{-1} &  \underline{ {\sigma_j}^{-1} {a_{ij}}^{-1} ({\sigma_j}^{-1}
{a_{ij}}^{-1} } {\sigma_j}^2 {a_{ij}} \sigma_j)   \\

 \ & \stackrel{M'_3}{=} & {a_{ij}} \underline{ {\sigma_j}^{-1} {a_{ij}}^{-1} {\sigma_j}^{-1}
{a_{ij}}^{-1} } {\sigma_j}^{-1} {\sigma_j}^2 {a_{ij}} \sigma_j    \\

 \ & \stackrel{M'_3}{=} & {a_{ij}} {a_{ij}}^{-1} {\sigma_j}^{-1}
{a_{ij}}^{-1} {\sigma_j}^{-1} \sigma_j {a_{ij}} \sigma_j    \\    [1.8mm]

 \ & = & 1.  \\
\end{array}$

\vspace{.1in}

\noindent Finally, for $P'_4$ and $r \in \{ 1, \ldots, m\}$ we have:

\vspace{.1in}

\noindent $\begin{array}{rcl}

 {a_{is}} & {a_{js}} & {a_{is}}^{-1} \ {a_{js}}^{-1} \ {a_{rs}} \ {a_{js}} \
 {a_{is}} \ {a_{js}}^{-1} \ {a_{is}}^{-1}    \\

 \ & \stackrel{M'_2}{=} &  (\sigma_{s-1} \ldots \sigma_j {a_{ij}}
{\sigma_j}^{-1} \ldots {\sigma_{s-1}}^{-1}) \cdot (\sigma_{s-1} \ldots \sigma_{j+1}
{\sigma_j}^2  {\sigma_{j+1}}^{-1} \ldots {\sigma_{s-1}}^{-1})  \\  [1.8mm]

\ & \ & \times \, (\sigma_{s-1} \ldots \sigma_j {a_{ij}}^{-1} {\sigma_j}^{-1} \ldots
{\sigma_{s-1}}^{-1}) \cdot
 (\sigma_{s-1} \ldots \sigma_{j+1} {\sigma_j}^{-2} {\sigma_{j+1}}^{-1} \ldots
{\sigma_{s-1}}^{-1})   \\   [1.8mm]

\ & \ & \times \, (\sigma_{s-1} \ldots \sigma_j {a_{rj}}
{\sigma_j}^{-1} \ldots {\sigma_{s-1}}^{-1}) \cdot (\sigma_{s-1} \ldots \sigma_{j+1}
{\sigma_j}^2  {\sigma_{j+1}}^{-1} \ldots {\sigma_{s-1}}^{-1})  \\  [1.8mm]

\ & \ & \times \, (\sigma_{s-1} \ldots \sigma_j {a_{ij}} {\sigma_j}^{-1} \ldots
{\sigma_{s-1}}^{-1}) \cdot
 (\sigma_{s-1} \ldots \sigma_{j+1} {\sigma_j}^{-2} {\sigma_{j+1}}^{-1} \ldots
{\sigma_{s-1}}^{-1})   \\   [1.8mm]

\ & \ & \times \, (\sigma_{s-1} \ldots \sigma_j {a_{ij}}^{-1} {\sigma_j}^{-1} \ldots
{\sigma_{s-1}}^{-1})   \\   [1.8mm]

 \ & = &  \sigma_{s-1} \ldots \sigma_{j+1} \underline{ \sigma_j {a_{ij}} {\sigma_j}^2
{a_{ij}}^{-1} {\sigma_j}^{-2} {a_{rj}} {\sigma_j}^2 {a_{ij}} {\sigma_j}^{-2} {a_{ij}}^{-1}
{\sigma_j}^{-1} } {\sigma_{j+1}}^{-1} \ldots {\sigma_{s-1}}^{-1}. \\
\end{array}$

\vspace{.1in}

\noindent On the other hand:

\vspace{.1in}

\noindent $\begin{array}{rcl}
 {a_{ij}}^{-1} &  \underline{ {a_{rs}} } & {a_{ij}} \ \stackrel{M'_2}{=}  \
\underline{ {a_{ij}}^{-1} \sigma_{s-1} \ldots \sigma_{j+1} } \sigma_j {a_{rj}}
{\sigma_j}^{-1} \underline{ {\sigma_{j+1}}^{-1} \ldots {\sigma_{s-1}}^{-1} {a_{ij}} }  \\

 \ & \stackrel{M'_1}{=} &
 \sigma_{s-1} \ldots \sigma_{j+1} \underline{ {a_{ij}}^{-1} \sigma_j {a_{rj}} {\sigma_j}^{-1}
 {a_{ij}} } {\sigma_{j+1}}^{-1}  \ldots {\sigma_{s-1}}^{-1}.   \\
\end{array}$

\vspace{.1in}

\noindent Again, it suffices to show that the underlined expressions in the last equation of
either side are equal. Indeed we have:

\vspace{.1in}

\noindent $\begin{array}{rcl}
\sigma_j & {a_{ij}} & {\sigma_j}^2 \, {a_{ij}}^{-1} \, {\sigma_j}^{-2} \, {a_{rj}} \, {\sigma_j}^2
 \, {a_{ij}} \, \underline{ {\sigma_j}^{-2} \, {a_{ij}}^{-1} \, {\sigma_j}^{-1} \, ({a_{ij}}^{-1} }
\, \sigma_j \, {a_{rj}}^{-1} \, {\sigma_j}^{-1} \, {a_{ij}})  \\

 \ & \stackrel{M'_3}{=} & \sigma_j {a_{ij}} {\sigma_j}^2 {a_{ij}}^{-1} {\sigma_j}^{-2}
{a_{rj}} {\sigma_j}^2 {a_{ij}} {a_{ij}}^{-1} {\sigma_j}^{-1} {a_{ij}}^{-1} {\sigma_j}^{-2}
\sigma_j {a_{rj}}^{-1} {\sigma_j}^{-1} {a_{ij}}  \\    [1.8mm]

 \ & = & \sigma_j {a_{ij}} {\sigma_j}^2 {a_{ij}}^{-1} {\sigma_j}^{-2} {a_{rj}}
\underline{ \sigma_j {a_{ij}}^{-1} {\sigma_j}^{-1} {a_{rj}}^{-1} } {\sigma_j}^{-1} {a_{ij}}  \\

 \ & \stackrel{(4)}{=} & \sigma_j {a_{ij}} {\sigma_j}^2 \underline{ {a_{ij}}^{-1} {\sigma_j}^{-1}
{a_{ij}}^{-1} {\sigma_j}^{-2} } {a_{ij}} \\

 \ & \stackrel{M'_3}{=} & \sigma_j {a_{ij}} {\sigma_j}^2 {\sigma_j}^{-2} {a_{ij}}^{-1}
{\sigma_j}^{-1} {a_{ij}}^{-1} {a_{ij}} \\    [1.8mm]

 \ & = & 1.      \\
\end{array}$

\vspace{.1in}

\noindent The case where $r \in \{ m+1, \ldots, m+n-1\}$ is completely analogous.  $\hfill \Box $

 \section{Irredundant presentation for $B_{m,n}$}

Looking now at the last presentation for $B_{m,n}$ we observe that  relations
$(2)$  may be seen as defining relations for $1 \leq i \leq m$ and $j \geq m+2$, namely:
$$\begin{array}{lcl}
a_{i,m+2}^{\pm} & := & \sigma_{m+1} a_{i,m+1}^{\pm} \sigma_{m+1}^{-1},  \\
a_{i,m+3}^{\pm} & := & \sigma_{m+2} \sigma_{m+1} a_{i,m+1}^{\pm} \sigma_{m+1}^{-1}
\sigma_{m+2}^{-1},  \\
     \ &  \vdots & \     \\
a_{i,m+n}^{\pm} & := & \sigma_{m+n-1} \ldots \sigma_{m+1} a_{i,m+1}^{\pm} \sigma_{m+1}^{-1} \ldots
\sigma_{m+n-1}^{-1}.   \\
\end{array}$$
 Therefore, we want to omit further these $a_{ij}$'s from the list of generators, and subsequently to
eliminate or simplify all relations involving these elements. Indeed, we have:

\begin{th}{ \ The following is a presentation for $B_{m,n}$:

\[ B_{m,n} = \left< \begin{array}{ll}  \begin{array}{l}
a_{1,m+1}, a_{2,m+1}, \ldots, a_{m,m+1},  \\
\sigma_{m+1}, \sigma_{m+2}, \ldots ,\sigma_{m+n-1} \\
\end{array} &
\left|
\begin{array}{l} \Sigma_1, \Sigma_2, (1'), (2'), (3'),       \\
                 \mbox{for all appropriate indices}   \\
\end{array}
\right.  \end{array} \right>,  \]

\noindent where:

\vspace{.1in}
\[\begin{array}{lrcl}
(1') \,  & {\sigma_k}^{-1}{a_{i,m+1}}^{\pm}{\sigma_k}   & = &  {a_{i,m+1}}^{\pm}, \ \ k \geq m+2,
\\
(2') \,  & {a_{i,m+1}}^{\pm} {\sigma_{m+1}} {a_{i,m+1}} {\sigma_{m+1}} & = &
 {\sigma_{m+1}} {a_{i,m+1}} {\sigma_{m+1}} {a_{i,m+1}}^{\pm},      \\

(3') \,  & {a_{i,m+1}}^{\pm} ({\sigma_{m+1}} {a_{r,m+1}}^{\pm} {\sigma_{m+1}}^{-1})  & = &
({\sigma_{m+1}} {a_{r,m+1}}^{\pm} {\sigma_{m+1}}^{-1}) {a_{i,m+1}}^{\pm}, \ \ r<i.  \\
 \end{array}\]
 } \end{th}

\bigbreak
\noindent {\bf Proof} \ Relations $(1')$ are a special case of relations $(1)$, relations $(2')$
are a special case of relations $(3)$ and relations $(3')$ are a special case of relations $(4)$.
Note that relations $(3')$ are equivalent to
$$ {a_{i,m+1}}^{\pm} {a_{r,m+2}}^{\pm}  =  {a_{r,m+2}}^{\pm} {a_{i,m+1}}^{\pm}, \ \ r<i. $$
We show the sufficiency of the new relations by  examining each case.

\vspace{.1in}

\noindent For relations $(1)$ and for $m+1 \leq k \leq j-2$ we have:

\vspace{.1in}

\noindent $\begin{array}{rcl}
 {\sigma_k}^{-1}  \underline{ {a_{ij}}^{\pm} }  {\sigma_k} & \stackrel{(2)}{=} &
\underline{ {\sigma_k}^{-1} (\sigma_{j-1} \ldots \sigma_{k+2} }\sigma_{k+1}
 \ldots \sigma_{m+1} {a_{i,m+1}}^{\pm}   \\   [1.8mm]

\ & \ & \times \, {\sigma_{m+1}}^{-1} \ldots {\sigma_{k+1}}^{-1} \underline{
 {\sigma_{k+2}}^{-1} \ldots {\sigma_{j-1}}^{-1}) {\sigma_k} }   \\

 \ &  \stackrel{\Sigma_1}{=} & (\sigma_{j-1} \ldots \sigma_{k+2}) \underline{ {\sigma_k}^{-1}
\sigma_{k+1} \sigma_k } \ldots \sigma_{m+1} {a_{i,m+1}}^{\pm}    \\   [1.8mm]

\ & \ & \times \, {\sigma_{m+1}}^{-1} \ldots
\underline{ {\sigma_k}^{-1} {\sigma_{k+1}}^{-1} \sigma_k } ({\sigma_{k+2}}^{-1} \ldots
{\sigma_{j-1}}^{-1})    \\

 \ &  \stackrel{\Sigma_2}{=} & (\sigma_{j-1} \ldots \sigma_{k+2}) \sigma_{k+1} \sigma_k
 \underline{  {\sigma_{k+1}}^{-1}  \sigma_{k-1} \ldots \sigma_{m+1} } {a_{i,m+1}}^{\pm}
   \\   [1.8mm]

\ & \ & \times \, \underline{ {\sigma_{m+1}}^{-1} \ldots {\sigma_{k-1}}^{-1} \sigma_{k+1} }
{\sigma_k}^{-1} {\sigma_{k+1}}^{-1} ({\sigma_{k+2}}^{-1} \ldots {\sigma_{j-1}}^{-1})    \\

 \ &  \stackrel{\Sigma_1}{=} & (\sigma_{j-1} \ldots \sigma_{m+1})
  \underline{ {\sigma_{k+1}}^{-1} {a_{i,m+1}}^{\pm} \sigma_{k+1} }
({\sigma_{m+1}}^{-1} \ldots {\sigma_{j-1}}^{-1})    \\

 \ &  \stackrel{(1')}{=} & \sigma_{j-1} \ldots \sigma_{m+1} {a_{i,m+1}}^{\pm} {\sigma_{m+1}}^{-1}
\ldots {\sigma_{j-1}}^{-1}, \ (k+1 \geq m+2)  \\

 \ &  \stackrel{(2)}{=} & {a_{ij}}^{\pm}.  \\
\end{array}$

\vspace{.1in}

\noindent Further, for relations $(1)$ and for $k \geq j+1$ we have:

\vspace{.1in}

\noindent $\begin{array}{rcl}
 {\sigma_k}^{-1}\underline{ {a_{ij}}^{\pm} } {\sigma_k}  &  \stackrel{(2)}{=} & \underline{
{\sigma_k}^{-1} (\sigma_{j-1} \ldots \sigma_{m+1} } {a_{i,m+1}}^{\pm} \underline{ {\sigma_{m+1}}^{-1}
\ldots {\sigma_{j-1}}^{-1}) {\sigma_k} }   \\

 \ &  \stackrel{\Sigma_1}{=} & (\sigma_{j-1} \ldots \sigma_{m+1}) \underline{ {\sigma_k}^{-1}
{a_{i,m+1}}^{\pm} {\sigma_k} } ({\sigma_{m+1}}^{-1} \ldots {\sigma_{j-1}}^{-1})   \\

 \ &  \stackrel{(1')}{=} & \sigma_{j-1} \ldots \sigma_{m+1} {a_{i,m+1}}^{\pm} {\sigma_{m+1}}^{-1}
\ldots {\sigma_{j-1}}^{-1})    \\

 \ &  \stackrel{(2)}{=} & {a_{ij}}^{\pm}.   \\
\end{array}$

\vspace{.1in}

\noindent In the last two equations we used that $k \geq m+2$. For relations $(3)$ we have:

\vspace{.1in}

\noindent $\begin{array}{rcl}
 {\sigma_j}\underline{ {a_{ij}} } {\sigma_j} \underline{ {a_{ij}}^{\pm} }  & \stackrel{(2)}{=} &
{\sigma_j} (\sigma_{j-1} \ldots \sigma_{m+1} {a_{i,m+1}}  \underline{ {\sigma_{m+1}}^{-1} \ldots
{\sigma_{j-1}}^{-1}) {\sigma_j} } \\   [1.8mm]

\ & \ & \times \, \underline{ ( \sigma_{j-1} \ldots \sigma_{m+1} } {a_{i,m+1}}^{\pm}
{\sigma_{m+1}}^{-1} \ldots {\sigma_{j-1}}^{-1} )    \\

 \ &  \stackrel{\Sigma_2}{=} & {\sigma_j} \ldots \sigma_{m+1} {a_{i,m+1}}
\underline{ {\sigma_j} } \ldots \underline{ \sigma_{m+2} } \sigma_{m+1} \underline{
{\sigma_{m+2}}^{-1} } \ldots \underline{ {\sigma_j}^{-1} }   \\   [1.8mm]

\ & \ & \times \,{a_{i,m+1}}^{\pm} {\sigma_{m+1}}^{-1} \ldots {\sigma_{j-1}}^{-1}    \\

 \ &  \stackrel{\Sigma_1, (1')}{=} & \underline{ {\sigma_j} ({\sigma_{j-1}}{\sigma_j}) \ldots
(\sigma_{m+1} \sigma_{m+2}) }{a_{i,m+1}} \sigma_{m+1} {a_{i,m+1}}^{\pm}  \\   [1.8mm]

\ & \ & \times \, ({\sigma_{m+2}}^{-1} {\sigma_{m+1}}^{-1}) \ldots ({\sigma_j}^{-1}
{\sigma_{j-1}}^{-1})      \\

 \ &  \stackrel{\Sigma_2}{=} & ({\sigma_{j-1}}{\sigma_j}) \ldots (\sigma_{m+1} \sigma_{m+2})
\underline{ \sigma_{m+1} {a_{i,m+1}} \sigma_{m+1} {a_{i,m+1}}^{\pm} }  \\   [1.8mm]

\ & \ & \times \, ({\sigma_{m+2}}^{-1} {\sigma_{m+1}}^{-1})  \ldots ({\sigma_j}^{-1}
{\sigma_{j-1}}^{-1})      \\

 \ &  \stackrel{(2')}{=} & ({\sigma_{j-1}}{\sigma_j}) \ldots (\sigma_{m+1} \sigma_{m+2})
{a_{i,m+1}}^{\pm} \sigma_{m+1} {a_{i,m+1}}  \\   [1.8mm]

\ & \ & \times \, \underline{ \sigma_{m+1} ({\sigma_{m+2}}^{-1}
{\sigma_{m+1}}^{-1}) \ldots ({\sigma_j}^{-1} {\sigma_{j-1}}^{-1}) }     \\

 \ &  \stackrel{\Sigma_2}{=} & ({\sigma_{j-1}} \underline{ {\sigma_j} }) \ldots (\sigma_{m+1}
\underline{ \sigma_{m+2} }) {a_{i,m+1}}^{\pm} \sigma_{m+1} {a_{i,m+1}}  \\   [1.8mm]

\ & \ & \times \, (\underline{ {\sigma_{m+2}}^{-1} } {\sigma_{m+1}}^{-1})
 \ldots (\underline{ {\sigma_j}^{-1} } {\sigma_{j-1}}^{-1}) \sigma_j     \\
\end{array}$

\vspace{.1in}

\noindent $\begin{array}{rcl}
 \ \ \ \ \  &  \stackrel{\Sigma_1, (1')}{=} & (\sigma_{j-1} \ldots \sigma_{m+1}) {a_{i,m+1}}^{\pm}
\underline{ {\sigma_j} \ldots \sigma_{m+2} \sigma_{m+1} {\sigma_{m+2}}^{-1} \ldots {\sigma_j}^{-1} }
  \\   [1.8mm]

\ & \ & \times \, {a_{i,m+1}} ({\sigma_{m+1}}^{-1} \ldots {\sigma_{j-1}}^{-1}) \sigma_j     \\

 \ &  \stackrel{\Sigma_2}{=} & \underline{ (\sigma_{j-1} \ldots \sigma_{m+1}) {a_{i,m+1}}^{\pm}
{\sigma_{m+1}}^{-1} \ldots {\sigma_{j-1}}^{-1} } \sigma_j \\   [1.8mm]

\ & \ & \times \,  \underline{ \sigma_{j-1} \ldots \sigma_{m+1} {a_{i,m+1}} } \underline{
({\sigma_{m+1}}^{-1} \ldots {\sigma_{j-1}}^{-1}) } \sigma_j     \\

 \ &  \stackrel{(2)}{=} &  {a_{ij}}^{\pm} \sigma_j {a_{ij}} \sigma_j.     \\
\end{array}$

\vspace{.1in}

\noindent Finally, for relations $(4)$ we have:

\vspace{.1in}

\noindent $\begin{array}{rcl}
\underline{ {a_{ij}}^{\pm} } & (\sigma_j & \underline{ {a_{rj}}^{\pm} } \ {\sigma_j}^{-1}) \
\stackrel{(2)}{=} \ (\sigma_{j-1} \ldots \sigma_{m+1} {a_{i,m+1}}^{\pm}
\underline{ {\sigma_{m+1}}^{-1} \ldots {\sigma_{j-1}}^{-1}) \sigma_j }    \\   [1.8mm]

\ & \ & \times \, \underline{ (\sigma_{j-1} \ldots \sigma_{m+1} } {a_{r,m+1}}^{\pm}
 {\sigma_{m+1}}^{-1} \ldots {\sigma_{j-1}}^{-1}) {\sigma_j}^{-1}    \\

 \ &  \stackrel{\Sigma_2}{=} & \sigma_{j-1} \ldots \sigma_{m+1} {a_{i,m+1}}^{\pm}
 \underline{ \sigma_j } \ldots \underline{ \sigma_{m+2} } \sigma_{m+1} \underline{
{\sigma_{m+2}}^{-1} } \ldots \underline{ {\sigma_j}^{-1} } {a_{r,m+1}}^{\pm}     \\   [1.8mm]

\ & \ & \times \, {\sigma_{m+1}}^{-1} \ldots {\sigma_j}^{-1}   \\

 \ &  \stackrel{\Sigma_1, (1')}{=} & ({\sigma_{j-1}}{\sigma_j}) \ldots (\sigma_{m+1} \sigma_{m+2})
{a_{i,m+1}}^{\pm} \sigma_{m+1} {a_{r,m+1}}^{\pm}  \\   [1.8mm]

\ & \ & \times \, \underline{ ({\sigma_{m+2}}^{-1}
{\sigma_{m+1}}^{-1}) \ldots ({\sigma_j}^{-1} {\sigma_{j-1}}^{-1}) {\sigma_j}^{-1} }    \\

\ &  \stackrel{\Sigma_2}{=} & ({\sigma_{j-1}}{\sigma_j}) \ldots (\sigma_{m+1} \sigma_{m+2})
\underline{ {a_{i,m+1}}^{\pm} (\sigma_{m+1} {a_{r,m+1}}^{\pm} {\sigma_{m+1}}^{-1}) }  \\   [1.8mm]

\ & \ & \times \, ({\sigma_{m+2}}^{-1} {\sigma_{m+1}}^{-1}) \ldots ({\sigma_j}^{-1}
{\sigma_{j-1}}^{-1})    \\

\ &  \stackrel{(3')}{=} & \underline{ ({\sigma_{j-1}}{\sigma_j}) \ldots (\sigma_{m+1} \sigma_{m+2})
\sigma_{m+1} } {a_{r,m+1}}^{\pm} {\sigma_{m+1}}^{-1} {a_{i,m+1}}^{\pm}   \\   [1.8mm]

\ & \ & \times \, ({\sigma_{m+2}}^{-1} {\sigma_{m+1}}^{-1}) \ldots ({\sigma_j}^{-1}
{\sigma_{j-1}}^{-1})    \\

\ &  \stackrel{\Sigma_2}{=} & {\sigma_j} ({\sigma_{j-1}} \underline{ {\sigma_j} }) \ldots
(\sigma_{m+1} \underline{ \sigma_{m+2} }) {a_{r,m+1}}^{\pm} {\sigma_{m+1}}^{-1} {a_{i,m+1}}^{\pm}
\\   [1.8mm]

\ & \ & \times \, (\underline{ {\sigma_{m+2}}^{-1} } {\sigma_{m+1}}^{-1}) \ldots (\underline{ {\sigma_j}^{-1} }
{\sigma_{j-1}}^{-1})    \\

\ &  \stackrel{\Sigma_1, (1')}{=} & {\sigma_j}  \ldots \sigma_{m+1} {a_{r,m+1}}^{\pm}
\underline{ {\sigma_j} \ldots \sigma_{m+2} {\sigma_{m+1}}^{-1} {\sigma_{m+2}}^{-1} \ldots
{\sigma_j}^{-1} }   \\   [1.8mm]

\ & \ & \times \, {a_{i,m+1}}^{\pm} {\sigma_{m+1}}^{-1} \ldots {\sigma_{j-1}}^{-1}    \\

\ &  \stackrel{\Sigma_2}{=} & {\sigma_j} \underline{ \sigma_{j-1} \ldots \sigma_{m+1}
{a_{r,m+1}}^{\pm} {\sigma_{m+1}}^{-1} \ldots {\sigma_{j-1}}^{-1} } {\sigma_j}^{-1}  \\   [1.8mm]

\ & \ & \times \, \underline{ \sigma_{j-1} \ldots \sigma_{m+1} {a_{i,m+1}}^{\pm} {\sigma_{m+1}}^{-1}
\ldots {\sigma_{j-1}}^{-1} }    \\

 \ &  \stackrel{(2)}{=} &  \sigma_j ({a_{rj}}^{\pm} {\sigma_j}^{-1} {a_{ij}}^{\pm} ). \hfill
\Box    \\
 \end{array}$

\bigbreak

\noindent The system of generators in the last presentation of $B_{m,n}$ is irredundant, in the sense
that no proper subset of it can generate $B_{m,n}$. In order now to simplify the notation
we will relabel the generators  $a_{1,m+1}, \ldots, a_{m,m+1}, \sigma_{m+1}, \ldots ,\sigma_{m+n-1}$
to  $a_1, \ldots, a_m, \sigma_1, \ldots ,\sigma_{n-1}$ accordingly, to obtain the
following, final presentation for $B_{m,n}$:

\[ B_{m,n} = \left< \begin{array}{ll}  \begin{array}{l}
a_1, \ldots, a_m,  \\
\sigma_1, \ldots ,\sigma_{n-1}  \\
\end{array} &
\left|
\begin{array}{l} \sigma_k \sigma_j=\sigma_j \sigma_k, \ \ |k-j|>1   \\
\sigma_k \sigma_{k+1} \sigma_k = \sigma_{k+1} \sigma_k \sigma_{k+1}, \ \  1 \leq k \leq n-1  \\
{a_i} \sigma_k = \sigma_k {a_i}, \ \ k \geq 2, \   1 \leq i \leq m,    \\
 {a_i} \sigma_1 {a_i} \sigma_1 = \sigma_1 {a_i} \sigma_1 {a_i}, \ \ 1 \leq i \leq m  \\
 {a_i} (\sigma_1 {a_r} {\sigma_1}^{-1}) =  (\sigma_1 {a_r} {\sigma_1}^{-1})  {a_i}, \ \ r < i.
\end{array} \right.  \end{array} \right>.  \]

\begin{rem}{\rm \ It is worth mentioning that the above presentation of $B_{m,n}$ is very similar to
that of the Artin braid group associated to the Dynkin diagram below.

\begin {center}
\leavevmode \epsfxsize =8cm \epsffile{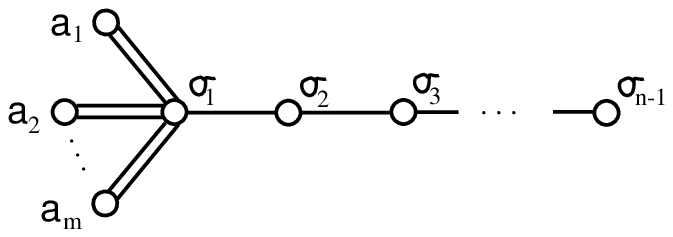}
\end {center}

\noindent In the diagram the single bonds mean relations of degree
$3$, the double bonds relations of degree $4$, and if two
generators are not connected by a bond they commute. The two
presentations differ only in the last relation, which in the case
of the Artin group (cf. [2]) is a mere commutation relation
between $a_i$ and $a_r$. Nevertheless, in the case $m=1$,
$B_{1,n}$ is the Artin group of type $\cal B$. } \end{rem}

\section{The cosets $C_{m,n}$}

In this section the word {\it knot} will be used to mean {\it knots and links}.

\vspace{.1in}
Let now $S^3 \backslash K$ be the  complement of the oriented knot
$K$ in $S^3$. Obviously, $S^3 \backslash K$ can be represented in
$S^3$ by the knot $K$. By classical results of Lickorish and
Wallace a closed, connected, orientable 3-manifold can be obtained
(not uniquely) by surgery along an integer-framed knot in $S^3$,
so it can be represented in $S^3$ by this knot. We shall denote by
$M$ either a knot complement or a c.c.o. 3-manifold. Then, by
fixing $M$ we may also fix a knot in $S^3$ representing $M$, and
this knot may be assumed to be a closed braid, say $\widehat{B}$.
It is shown in [7] that knots in these spaces can be represented
by `mixed braids', which contain the braid $B$ as a fixed
subbraid. More precisely, we have the following:

\begin{defn}{ \ A mixed braid is a special element of $B_{m+n}$ consisting of two disjoint orbits of
strands, such that the subbraid forming the one orbit consists of the first $m$ strands and it is a
fixed element of $B_m$}. \end{defn}

If now  the manifold $M$  is a handlebody of genus $m$, knots in $M$ can clearly be  represented by
elements of $B_{m,n}$. Further, if  $M$ is the complement of the $m$-unlink or a connected sum of $m$
lens spaces of  type $L(p,1)$, then knots in $M$ are represented  by   elements  of $B_{m,n}$,
for $n \in \NN$. In the special case where $M$ is the complement of the trivial knot or
$L(p,1)$, knots in $M$ are represented by  elements of $B_{1,n}$, the Artin group of type $\cal B$
 (cf. [4, 5, 6]).  These  are rather special cases of
 knot complements or c.c.o. 3-manifolds.

\bigbreak

In the generic case the subbraid representing $M$ will not be the identity braid (for an example see
figure 3(a)). In the generic case the multiplication of two mixed braids in $B_{m+n}$ related to $M$
is not a closed operation, since concatenation does not preserve the structure of the manifold.

\begin {center}
\leavevmode \epsfxsize =13cm \epsffile{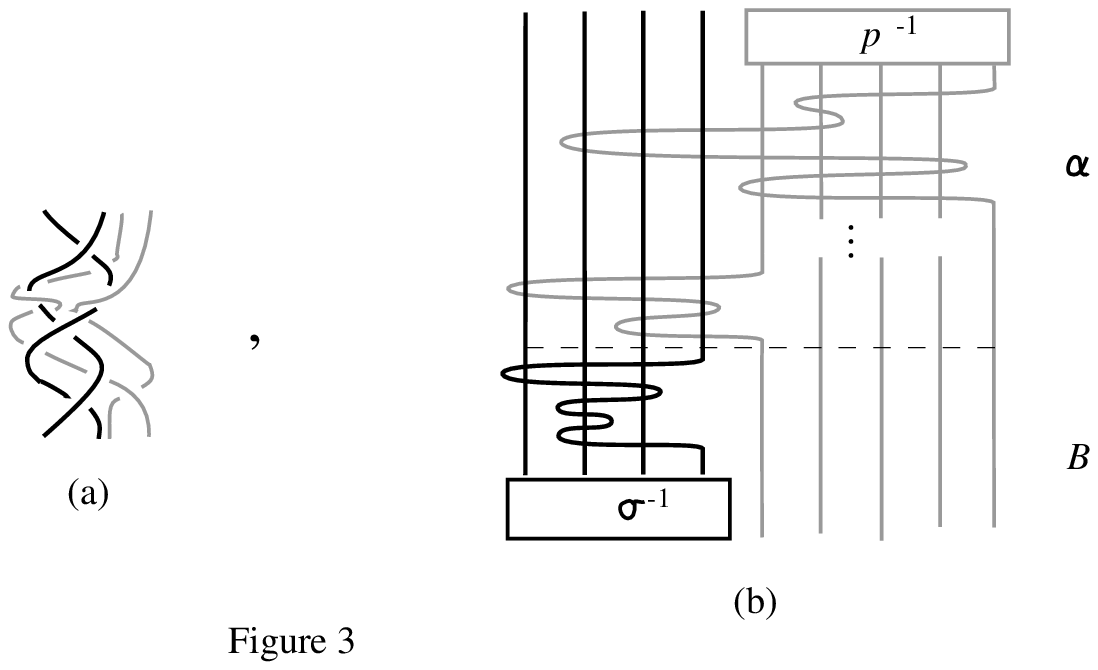}
\end {center}

\noindent The following proposition shows that, nevertheless, we still have
braid structures in $M$. Indeed, let  $\bigcup_{n=1}^{\infty} C_{m,n}$  denote the disjoint union of
the cosets of all mixed braids associated to a generic $M$. Then we have the following:

\bigbreak

\begin{prop}{ \ \  For a fixed $M$, $C_{m,n}$ is a  coset of $B_{m,n}$ in $B_{m+n}$.}\end{prop}

\bigbreak
\noindent {\bf Proof} \ Let  $A\in C_{m,n}$. We shall show that $A$ can
be written as a product  $\alpha\cdot B$, where  $\alpha\in
B_{m,n}$  and  $B$ is the fixed braid representing $M$ in $S^3$. Indeed, we notice first that
by symmetry, Artin's combing for pure braids can be also applied starting from the last strand of a
pure braid. So, we multiply $A$ from the top with a braid $p$ on the last $n$ strands and with a
braid $\sigma$ on the first $m$ strands, such that $pA\sigma$ is a pure braid in  $B_{m+n}$. Then
we apply to it Artin's canonical form from the end. This will separate $A$ into two parts, one
being an element of  $B_{m,n}$ and the other being the fixed braid $B$ embedded in $B_{m+n}$ (see
figure 3(b)).
 $\hfill \Box  $

\bigbreak

A final comment is now due.

\begin{rem} {\rm \ The groups $B_{m,n}$ and their cosets $C_{m,n}$ will be used for yielding an
algebraic version of Markov's theorem for isotopy of knots in knot complements and 3-manifolds. For
the purpose of constructing knot invariants following the line of  Jones-Ocneanu one can use the
irredundant presentation of $B_{m,n}$ for considering appropriate quotient algebras that satisfy a
quadratic skein relation for the $\sigma_i$'s.  } \end{rem}


\bigbreak
\noindent {\sc S.L.: Bunsenstrasse 3-5, Mathematisches Institut, G\"{o}ttingen Universit\"{a}t,
37073 G\"{o}ttingen, Germany. E-mail: sofia@cfgauss.uni-math.gwdg.de }


\begin{thebibliography}{131}

\bibitem{B} {J.S. Birman}, ``Braids, links and mapping class groups'',  Ann. of Math. Stud.
            {\bf 82}, Princeton University Press, Princeton, 1974.

\bibitem{BS} {E. Brieskorn, K. Saito}, Artin-Gruppen und Coxeter-Gruppen,
              {\it Inventiones Math.} {\bf 17}, 245--271 (1972).

\bibitem{Jo}  {D.L. Johnson}, ``Presentations of Groups'',  LMS Student Texts {\bf 15}, 1990.

\bibitem{La} {S. Lambropoulou}, ``A study of braids in 3--manifolds'',
             Ph.D. thesis, Warwick, 1993.

\bibitem{La1} {S. Lambropoulou},  Solid torus links and Hecke algebras of B-type,
             Proceedings of the Conference on Quantum Topology, D.N. Yetter ed., {\it World
             Scientific Press}, 1994.

\bibitem{La2} {S. Lambropoulou},  Knot theory related to generalized and cyclotomic Hecke algebras
             of type {\cal B}, {\it J. Knot Theory and its Ramifications} {\bf Vol. 8, No. 5},
             621--658 (1999).

\bibitem{LR} {S.~Lambropoulou, C.P.~Rourke}, Markov's theorem in 3-manifolds, {\it Topology and
its Applications} {\bf 78}, 95--122 (1997).

\bibitem{S}  {A.B. Sossinsky},
               Preparation theorems for isotopy invariants of links in 3-manifolds,
               Quantum Groups, Proceedings, {\it Lecture Notes in Math.} {\bf 1510}, Springer-Verlag
               Berlin a.o. 354--362 (1992).

\bibitem{V}  {V.V. Vershinin},
               Homology of braid groups in handlebodies, Preprint No 96/06-2, Universit\'{e} de
             Nantes (1996).

\end{thebibliography}
\end{document}